\documentclass[11pt]{article}
\usepackage[letterpaper,twoside,outer=1.3in,vmargin=1.3in,]{geometry}
\usepackage{amsmath,amssymb,amsthm,amsfonts,enumerate,times,epsfig,color,xcolor, tikz}
\usepackage{thmtools}
\usepackage{thm-restate}
\usepackage{hyperref,cleveref,soul}
\usepackage[numbers]{natbib}

\usetikzlibrary{graphs, shapes, calc}
\usepackage{subcaption}  

\allowdisplaybreaks

\newcommand{\ignore}[1]{}

\DeclareMathOperator{\lin}{lin}

\DeclareMathOperator{\pete}{\mathbb{P}}

\newcommand{\1}{\mathbf{1}}
\newcommand{\0}{\mathbf{0}}

\newcommand{\cR}{\mathbb R}
\newcommand{\cZ}{\mathbb Z}

\newtheorem{theorem}{Theorem}[section]

\newtheorem{LE}[theorem]{Lemma}

\newtheoremstyle{plain-upright}
  {\topsep}   
  {\topsep}   
  {\upshape}  
  {}          
  {\bfseries} 
  {.}         
  { }         
  {}          

\theoremstyle{plain-upright}
\newtheorem{custom}[theorem]{}  

\newcounter{claim_nb}[theorem]
\newtheorem*{claim*}{Claim}
\newtheorem*{subclaim*}{Subclaim}

\newcounter{claim_nbs}[section]
\newcounter{subclaim_nb}[claim_nbs]
\newif\ifnotes\notesfalse

\ifnotes
\usepackage{color}

\renewcommand{\b}[1]{{\color{blue} #1}}
\renewcommand{\r}[1]{{\color{red} #1}}
\newcommand{\rt}[1]{{\color{red}{\st{#1}}}} 
\newcommand{\rb}[2]{\rt{#1}\b{#2}}

\else

\renewcommand{\b}[1]{#1}
\renewcommand{\r}[1]{}
\newcommand{\rt}[1]{}
\newcommand{\rb}[2]{#2}

\fi

\title{Integral bases, perfect matchings, and the Petersen graph}

\author{Ahmad Abdi \and Olha Silina}

\begin{document}

\maketitle

\begin{center}
\textit{Dedicated to the memory of Professor Murty}
\end{center}

\begin{abstract}
Let $G=(V,E)$ be a matching-covered graph, denote by $P$ its perfect matching polytope, and by $L$ the integer lattice generated by the integral points in $P$. In this paper, we give short, polyhedral proofs for two difficult results established by Lov\'{a}sz (1987), and by Carvalho, Lucchesi, and Murty (2002) in a series of three papers totaling over 120 pages. More specifically, we prove that $L$ has a lattice basis consisting solely of incidence vectors of some perfect matchings of~$G$, $2x\in L$ for all $x\in \lin(P)\cap \cZ^E$, and if $G$ has no Petersen brick then $L = \lin(P)\cap \cZ^E$. Our proof avoids major technical aspects of the previous proofs, the most important of these being a characterization of the dual lattice, and a `Petersen-brick-sensitive' ear decomposition result for matching-covered graphs. This is achieved by a novel study of the facial structure of the polytope $P$ and its relationship with the lattice $L$. It is also based on a first-of-its-kind polyhedral characterization of the Petersen graph.\\

\end{abstract}


\section{Introduction}\label{sec:intro}

Perfect matchings lie at the heart of combinatorial optimization and structural graph theory. Consider a \emph{matching-covered} graph $G=(V,E)$, that is, a graph where every edge appears in a perfect matching. Denote by $P(G)\subseteq \cR^E$ the polytope whose vertices are the incidence vectors of the perfect matchings of $G$, and by $L(G)\subseteq \cZ^E$ the (integer) lattice generated by the integral points in $P(G)$, i.e., $L(G)$ is the set of all integer linear combinations of $P(G)\cap \{0,1\}^E$. Understanding the relationship between $P(G)$, the combinatorial structure of $G$, and the algebra of $L(G)$ lies at the center of matching theory.

A seminal result in combinatorial optimization due to Edmonds~\cite{Edmonds65} is that $P(G)$ can be described by non-negativity inequalities: $x_e\geq 0,e\in E$, degree equations: $x(\delta(v))=1,v\in V$, and \emph{odd cut} inequalities amounting to $x(\delta(U))\geq 1$ for all odd-sized subsets $U\subset V$. Here, $x(A)$ denotes $\sum_{a\in A} x_a$ for an edge set $A\subseteq E$. Later Seymour~\cite{Seymour79} gave a graph-theoretic proof of this result using Tutte's characterization of perfect matchings~\cite{Tutte47}. In the same paper, Seymour connected the polyhedral study of $P(G)$ with $L(G)$. Through a long and technical argument, he proved that in a bridgeless cubic graph, the all-$2$s vector belongs to $L(G)$, and in fact the all-$1$s vector belongs to the lattice if $G$ has no Petersen minor. This, in turn, established weaker variants of the Berge-Fulkerson conjecture and Tutte's $4$-flow conjecture, respectively~\cite{Fulkerson71,Tutte66}. 

Many subsequent works have led to a deeper understanding of the structure of $P(G)$ and $L(G)$, as well as various connections between these objects.
Naddef~\cite{Naddef82} and Edmonds, Pulleyblank, and Lov\'{a}sz~\cite{Edmonds82} determined the dimension of $P(G)$ based on structural properties of $G$, showing how combinatorial decomposition manifests itself in the geometry of the polytope. These and related lines of questions led to the development of a fascinating area known as \emph{matching theory}~\cite{LP09,LM24}. In a seminal paper of the area, Lov\'{a}sz \cite{Lovasz87} gave a far-reaching extension of Seymour's result and characterized the lattice $L(G)$ itself. He proved that $2x\in L(G)$ for all $x\in \lin(P(G))\cap \cZ^E$, and if $G$ has no Petersen brick (defined later) then $L(G) = \lin(P(G))\cap \cZ^E$. Here, $\lin(S)$ denotes the linear hull of $S$ for $S\subseteq \cR^E$. The core of his argument was proving a difficult lemma that characterized the dual lattice for a non-Petersen brick. This lengthy proof relied on key tools from matching theory, including the notions of barrier cuts, 2-separation cuts, and ear decompositions. 

Over a decade later, Carvalho, Lucchesi, and Murty proved that $L(G)$ has a \emph{lattice basis} $B$ consisting solely of incidence vectors of some perfect matchings of~$G$. That is, $B$ is a (linear) basis for $\lin(L(G))$, and every vector in $L(G)$ is an integer linear combination of $B$. This fact was featured as a key application in a series of three excellent papers~\cite{CLM1,CLM2,CLM3}, totaling over 120 pages. There, the authors develop a special type of ear decomposition of the graph $G$ that, vaguely speaking, is sensitive to the number of bricks and Petersen bricks. This in turn gave an affirmative answer to a question previously raised by Murty~(\cite{Murty94}, Problem 7.3) about the existence of such a basis.

\subsection{The Main Theorem and the Integral Basis Theorem}

In this paper, we give relatively short, polyhedral proofs of the theorems mentioned above due to Lov\'{a}sz, and Carvalho, Lucchesi, and Murty. To formally state these theorems, it is necessary to introduce some notation and concepts.

Let $G=(V,E)$ be a matching-covered graph. For an odd cut $C$ and an edge $e\in E$, denote by $P(G;C)$ and $P_e(G)$ the faces $P(G)\cap \{x:x(C)=1\}$ and $P(G)\cap \{x:x_e=0\}$ of $P(G)$, respectively. Note that $P_e(G)\subsetneq P(G)$ as $G$ is matching-covered.

\paragraph{Separating cuts.}
Let $C=\delta(X)$ be an odd cut. We say that $C$ is \emph{separating}, or \emph{contractible}, in $G$ if $1<|X|<|V|-1$, and each $e\in E$ belongs to some perfect matching $M$ of $G$ with $|M\cap C|=1$. The second condition is equivalent to the polyhedral condition that $P(G;C)\not\subseteq P_e(G)$ for any edge $e\in E$. 
Given a separating cut $C=\delta(X)$, we refer to $G/X,G/\bar{X}$ as the \emph{$C$-contractions} or \emph{cut-contractions} of $G$ \b{where $\bar{X}$ denotes the complement $V\setminus X$}. We denote the \emph{contraction vertices} in a cut-contraction $G/X$ by the corresponding lower case letter $x$. 
An odd cut $C$ is separating if, and only if, both $C$-contractions are matching-covered. Observe further that if $\delta(X),\delta(Y)$ are odd cuts of $G$ such that $X\subset Y$, where $\delta(Y)$ is separating in $G$, and $\delta(X)$ is separating in $G/\bar{Y}$, then $\delta(X)$ is separating in $G$. This is because any perfect matching of a separating cut-contraction can be extended to a perfect matching of the original graph.

\begin{figure}[h!]\label{pic:petersen}
\centering
\begin{tikzpicture}[scale=1.4, every node/.style={circle, draw, fill=black, inner sep=2pt}, spoke/.style={line width=5pt, color=yellow!80!black, opacity=0.5}]

\def\R{2}       
\def\r{1}       
\def\shift{18}

\foreach \i in {1,...,5} {
  \node (o\i) at ({72*(\i-1)+\shift}:\R) {};
}

\foreach \i in {1,...,5} {
  \node (i\i) at ({72*(\i-1)+\shift}:\r) {};
}

\foreach \i [evaluate=\i as \j using {mod(\i,5)+1}] in {1,...,5} {
  \draw (o\i) -- (o\j);
}

\foreach \i/\j in {1/3, 3/5, 5/2, 2/4, 4/1} {
  \draw (i\i) -- (i\j);
}

\foreach \i in {1,...,5} {
  \draw (o\i) -- (i\i);  
}

\foreach \i in {1,...,5} {
  \draw[spoke] (o\i) -- (i\i);
}
\end{tikzpicture}
\caption{The Petersen graph. Any cut separating two $5$-cycles (e.g., the one highlighted) is a separating facet-defining cut.}
\label{fig:petersen}
\end{figure}
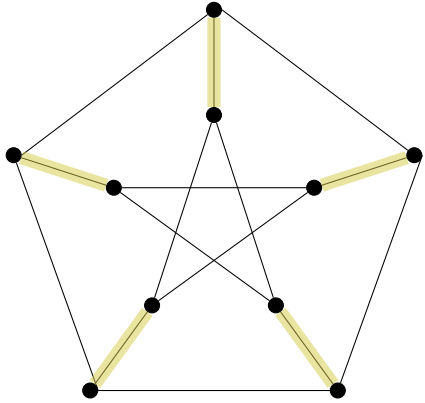

\paragraph{Bricks.}
A \emph{tight cut} is an odd cut $C=\delta(X)$ such that $1<|X|<|V|-1$, and every perfect matching of $G$ intersects $C$ exactly once. The second condition is equivalent to the polyhedral condition $P(G;C)=P(G)$. Tight cuts are clearly separating. A matching-covered graph with no tight cut is a \emph{brick} if it is non-bipartite, and is a \emph{brace} otherwise. 
Cut-contractions along tight cuts repeatedly give rise to a binary tree rooted at $G$ whose leaves correspond to a \emph{tight cut decomposition} of $G$ into bricks and braces. In fact, the list of bricks and braces in a tight cut decomposition is unique up to the multiplicity of the edges in each brick and brace~\cite{Lovasz87}. 
One specific brick plays an important role in matching theory, namely, the Petersen brick.
It can be checked that the Petersen graph (see \Cref{fig:petersen} for an illustration) is a brick. Throughout the work, we allow graphs to have parallel edges, but not loops; with this in mind, we refer to any graph whose underlying simple graph is the Petersen graph, as a \emph{Petersen brick}. 

\medskip
The main result of this work is a polyhedral proof of the following theorem.

\begin{theorem}[\cite{CLM3,Lovasz87}]\label{lattice-basis}
Let $G=(V,E)$ be a matching-covered graph, let $L:=L(G)$ and $\bar{L}:=\lin(P(G))\cap \cZ^E$. Then $L$ has a lattice basis consisting solely of perfect matchings of $G$. Furthermore, if $G$ has $p$ Petersen bricks in its tight cut decomposition, then 
$$L = \bar{L} \cap \left\{x : x(A_i) \equiv 0 \pmod{2},~\forall i\in [p]\right\},$$
where each $A_i$ is the edge set of some $5$-cycle of the $i\textsuperscript{th}$ Petersen brick. In particular, if $p=0$ then $L=\bar{L}$, and if $p\geq 1$ then $2x\in L$ for all $x\in \bar{L}$.
\end{theorem}

Above, $[p]$ denotes $\{1,2,\ldots,p\}$ for $p\geq 1$, and $\emptyset$ for $p=0$. Similar to the approach of \cite{CLM3}, we make extensive use of tight and more generally separating cuts. However, apart from these, our proof mostly circumvents the lengthy graph-theoretic arguments, and in particular has no dependence on matching theoretic notions such as the uniqueness of the bricks and braces, $2$-separation cuts, braces, removable edges, and ear decompositions; our proof has minimal dependence on the notions of barriers and bricks. It also eliminates the need to study the dual lattice through matching-integral vectors. Instead, our proof is based on polyhedral, or otherwise polyhedrally-driven graph-theoretic arguments that study the relationship between the lattice and the facial structure of the polytope $P(G)$. 

\paragraph{The Integral Basis Theorem.} A matching-covered graph is \emph{Petersen-free} if it has no Petersen brick in its tight cut decomposition. A key notion needed for our proof is that of an \emph{integral basis} for a rational linear subspace, which is a basis $B$ \b{consisting of integer vectors} such that every integral vector in the subspace is an integer linear combination of $B$. The crux of \Cref{lattice-basis} is the following special case. 

\begin{theorem}\label{integral-basis}
Let $G=(V,E)$ be a Petersen-free matching-covered graph. Then $\lin(P(G))$ has an integral basis consisting solely of incidence vectors of perfect matchings. In particular, $L(G)=\lin(P(G))\cap \cZ^E$.
\end{theorem}

\subsection{Polyhedral theory and the perfect matching polytope}

Before explaining the key ideas behind our proof, we first recall some concepts from polyhedral theory, and make connections to the perfect matching polytope. We also refer the reader to \cite{conforti2014integer} for more polyhedral theory. 

\paragraph{Primer on polyhedral theory.} Given a polyhedron $Q=\{x\in \cR^n: Ax\leq b, Cx = d\}$, a \emph{face} of $Q$ is any nonempty polyhedron obtained by intersecting $Q$ with $\{x:a^\top x= \beta\}$ for some inequality $a^\top x\leq \beta$ valid for $Q$. Such an inequality $a^\top x\leq \beta$ \emph{defines}, or \emph{exposes} the face $Q\cap \{x\in \cR^n: a^\top x=\beta\}$. 
We further classify faces based on a measure called \emph{dimension}: the (affine) dimension of a polyhedron $Q$ is the maximum number of affinely independent points in $Q$ minus one, denoted by $\dim(Q)$. 
\b{We will mostly be concerned with polytopes whose affine hull does not contain ${\bf 0}$; for such polytopes the dimension is the maximum number of \emph{linearly} independent points minus one.} 
A \emph{facet} of $Q$ is a face whose dimension is exactly $\dim(Q)-1$; furthermore, an inequality $a^\top x\leq \beta$ of $Ax\leq b$ is \emph{facet-defining} if it exposes a facet of $Q$. Facets play an essential role in the description of the polyhedron, as any polyhedron is described by its implicit equality constraints (which include $Cx=d$ and possibly a subset of $Ax\leq b$) and its facet-defining inequality constraints (which form a subset of $Ax\leq b$).

\paragraph{Dimension of the perfect matching polytope.}
Let $G=(V,E)$ be a matching-covered graph. The matching polytope $P(G)$ is known to have dimension $|E|-|V|+1-b(G)$, where $b(G)$ is the number of bricks of $G$ in a tight cut decomposition~\cite{Naddef82,Edmonds82}. Further, $b(G)=0$ if and only if $G$ is bipartite~\cite{Naddef82}. 

\paragraph{BvN graphs.} 
Recall that the perfect matching polytope $P(G)$ is described by three types of constraints: degree constraints, non-negativity constraints, and odd cut constraints. Since the degree constraints are equalities, there are only two types of facets in the description of $P(G)$, exposed by (a) $x_e\geq 0$ for some $e\in E$, or (b) $x(C)\geq 1$ for some odd cut $C\subseteq E$. Note that in general, the facets of $P(G)$ might be exposed by both types of inequalities at once. In this work, we distinguish a class of graphs where each facet is exposed by some non-negativity constraint. Formally, a graph $G=(V,E)$ is \emph{Birkhoﬀ–von Neumann (BvN)} if $P(G) = \left\{x\in \cR_{\geq 0}^E : x(\delta(v))=1,\,\forall v\in V\right\}$. This terminology was coined recently in \cite{Carvalho20}, but the class was studied and characterized \rb{much}{} earlier, see \cite{Carvalho04}. It is well-known that every bipartite graph is BvN, though some non-bipartite graphs such as $K_4$ can also be BvN.

\paragraph{Separating facet-defining cuts and non-BvN near-bricks.}
\rb{Suppose a matching-covered graph $G$ is \emph{not} BvN. By definition, this happens if, and only if, $P(G)$ has a facet that is not exposed by any non-negativity inequality $x_e\geq 0$. Such a facet must be}
{Suppose a matching-covered graph $G$ is \emph{not} BvN. By definition, this happens if and only if one of the following holds: \begin{itemize}
    \item[(i)] there is a tight cut $C$ such that $x(C)-1$ cannot be written as a linear combination of the trivial cuts $x(\delta(v))-1$, or
    \item[(ii)] $P(G)$ has a facet that is not exposed by any non-negativity inequality $x_e\geq 0$. 
\end{itemize}

If we further assume that $G$ is a \emph{near-brick}, i.e., has exactly one brick in its tight cut decomposition, then we must be in the second case as near-bricks have no tight cut satisfying (i). A facet satisfying this case must be} exposed by (a special type of) an odd cut inequality, denote it $x(C)\geq 1$. Observe that $C$ must be separating since otherwise, there would exist an edge $e\in E$ such that $P(G)\cap \{x:x(C)=1\}\subseteq P(G)\cap \{x:x_e= 0\}\subsetneq P(G)$ --- recall that $G$ is matching-covered. However, as $x(C)\geq 1$ exposes a facet of $P(G)$, the inequality $x_e\geq 0$ must expose the same facet, contradicting our choice of the facet.
Hence, such a cut $C$ must be separating; we refer to $C$ as a \emph{separating facet-defining cut} of $G$. Consequently, \b{a near-brick} $G$ is not BvN if, and only if, there exists a separating facet-defining cut. 

\subsection{Proof ingredients: the Petersen Graph Lemma and the Intersection Theorem}

Our proof of \Cref{integral-basis} is based on two key ingredients.

\paragraph{The Petersen Graph Lemma.} 
The following lemma provides a novel polyhedral characterization of the Petersen graph, distinguishing it from all other bricks.

\begin{LE}\label{petersen-LE}
Let $G=(V,E)$ be a brick, and $d$ the dimension of $P(G)$. Suppose every face of dimension $d-2$ is exposed by a non-negativity inequality. Then $|V|\leq 10$. Furthermore, $G$ is the Petersen graph, or $G$ is BvN, or there exist a perfect matching $M$ and a separating facet-defining cut $C$ such that $|M\cap C|=3$.
\end{LE}

\paragraph{The Intersection Theorem.} We shall also need the following theorem.

\begin{theorem}\label{intersection-theorem}
Let $G=(V,E)$ be a Petersen-free non-BvN near-brick. Then there exists a perfect matching $M$ and a separating facet-defining cut $C$ such that $|M\cap C|=3$.
\end{theorem}

This is a special case of a more general unpublished result by Campos and Lucchesi~\cite{Campos00}, claiming in particular that in $G$ as above, for every non-tight separating cut $C$, there is a perfect matching $M$ such that $|M\cap C|=3$. The Intersection Theorem is a slight variant of the main theorem in the two papers \cite{CLM1,CLM2}, which establishes in an arbitrary matching-covered graph the minimum intersection size between a perfect matching and a separating cut with intersection size at least three.

\paragraph{Proof outline.} 
The Main Theorem follows rather routinely from the Integral Basis Theorem. The proof of the latter proceeds by reducing the problem to bricks. If $G$ is BvN, then $P(G)$ has the `integer decomposition property', which we use to settle the claim. Otherwise, by the Intersection Theorem, there is a separating facet-defining cut $C$ and a perfect matching $M$ such that $|C\cap M|=3$. After potentially tweaking $(C,M)$, we then decompose $G$ into two smaller Petersen-free matching-covered graphs $G_1$ and $G_2$, find integral bases $B_1$ and $B_2$ for each and compose them, and then add the incidence vector of $M$ to obtain an integral basis for $G$.

The proof of the Intersection Theorem proceeds by first reducing to bricks. The key idea is to study the $(d-1)$- and $(d-2)$-dimensional faces of the perfect matching polytope, where $d=\dim(P(G))$. If $P(G)$ has a facet of the form $x_e\geq 0$, we then apply induction to $G\setminus e$ for a suitably chosen such edge $e$. If $P(G)$ has a $(d-2)$-dimensional face $F$ \emph{not} exposed by any $x_f\geq 0$, then roughly speaking, we apply induction to $G/U$ where $x(\delta(U))\geq 1$ is a facet `derived' from $F$. We then resort to the Petersen Graph Lemma to finish the proof.

\paragraph{Paper outline.} 
In \S\ref{sec:petersen-graph-LE}, we give a short proof of the Petersen Graph Lemma. In \S\ref{sec:ingredients-2}, we present further preliminary results, and in \S\ref{sec:intersection-theorem}, we prove the Intersection Theorem. Finally, in \S\ref{sec:integral-basis-theorem}, we prove the Integral Basis Theorem, and then obtain the Main Theorem as a consequence.

\section{Proof of the Petersen Graph Lemma}\label{sec:petersen-graph-LE}

Let $v:=|V|$ and $e:=|E|$, let $f$ be the number of facets and $t$ the number of \emph{$(d-2)$-faces}, i.e., faces of dimension $d-2$, of $P(G)$. As $G$ is a brick, it follows that $d = e-v$.
Observe that every facet of the $d$-dimensional polytope $P(G)$ has at least $d$ distinct adjacent facets, and each adjacency defines a distinct $(d-2)$-face. Furthermore, every $(d-2)$-face belongs to exactly $2$ facets. Thus, $
t \geq \frac{f \cdot d}{2}$. Given that every $(d-2)$-face is exposed by $x_e\geq 0$ for some edge $e\in E$, and each such edge defines at most one $(d-2)$-face, it follows that $e\geq t$. Furthermore, $f\geq d+1$ as $P(G)$ is $d$-dimensional, so $$
e\geq t\geq \frac{fd}{2}\geq {d+1 \choose 2} = {e-v+1 \choose 2}.
$$ 

Subsequently, $e\geq {e-v+1 \choose 2}$, which can be rewritten as $e+v\geq (e-v)^2$. Given that $G$ is a brick, we may assume every vertex $u$ has degree at least $3$ (otherwise, $u$ would have $1$ or $2$ neighbors, and as $G$ is a brick, $v\in \{2,4\}$). Subsequently, $2e-3v\geq 0$. Let us now solve the following convex minimization problem:
$$\min\left\{(x-y)^2 -x-y:-2x+3y\leq 0,10-y\leq 0\right\}.$$ We see that at the minimum, $x=15, y=10$, and the Lagrange multipliers for the inequalities are $\lambda=9/2$ and $\mu=5/2$, respectively, and the optimal value is $0$. Subsequently, the inequalities $e+v\geq (e-v)^2$ and $2e-3v\geq 0$ imply that $v\leq 10$, and if $v=10$ then $e=15$.

Suppose $G$ is not BvN, and there do not exist a perfect matching $M$ and a separating facet-defining cut $C$ such that $|M\cap C|=3$. We shall prove that $G$ is the Petersen graph.

\rb{We claim that $v=10$. If not, then $v\leq 8$. As $G$ is not BvN, $P(G)$ has a facet not exposed by a non-negativity inequality, say it is exposed by an odd cut $C=\delta(X)$. Then $C$ must be separating, so $3\leq |X|\leq |\bar{X}|\leq v-3$. Thus, $v\in \{6,8\}$ and $|X|=3$. Let $M$ be a perfect matching such that $|M\cap C|>1$. As $|X|=3$, it follows that $|M\cap C|=3$, a contradiction.}{
As $G$ is a non-BvN brick, $P(G)$ has a facet not exposed by a non-negativity inequality, say it is exposed by a separating facet-defining cut $C= \delta(X)$.
Since $C$ is not tight, there exists a perfect matching $M$ such that $|M\cap C|> 1$, hence $|M\cap C|\geq 5$. Since $v\leq 10$, this yields $|M|= |X|= |V \setminus X|= 5$ and $v= 10$.
}
Thus, \rt{$v=10$, and so} $e=15$ and $e+v=(e-v)^2$, implying that $e=t = \frac{fd}{2} = {d+1 \choose 2}$, and so $f=d+1=6$. In particular, $G$ is a cubic graph as $2e=3v$\rt{,and no facet is exposed by a non-negativity inequality as $e=t$}. \rt{Let $C=\delta(X)$ be one of the $6$ facet-defining cuts, which must be separating. As there is no perfect matching intersecting $C$ three times, it follows that $|X|=|\bar{X}|=5$, and $|C|\geq 5$.} As connected subgraphs, $G[X],G[\bar{X}]$ each contains at least $4$ edges, so $|C|\in \{5,7\}$ and both of $G[X],G[\bar{X}]$ are either $5$-cycles or $4$-paths. In the case of the latter, the $C$-contractions are not matching-covered, which is a contradiction as $C$ is separating. Thus, $|C|=5$ and both of $G[X],G[\bar{X}]$ are $5$-cycles. 
As there is no perfect matching intersecting $C$ three times, there is no $4$-cycle intersecting $C$ twice, implying in turn that $G$ is the Petersen graph, as required.\qed

\section{Ingredients for the Intersection Theorem}\label{sec:ingredients-2}

Let $G=(V,E)$ be a matching-covered graph, and let $d:=\dim(P(G))$. Let $C=\delta(X)$ be a separating cut, and let $G_1:=G/X$ and $G_2:=G/\bar{X}$. Suppose 
$$
B_1:=\{x^1,\ldots,x^{\hat{d}_1}\}\subseteq P(G_1)\cap \{0,1\}^{E(G_1)}
\quad\text{and}\quad
B_2:=\{y^1,\ldots,y^{\hat{d}_2}\}\subseteq P(G_2)\cap \{0,1\}^{E(G_2)}
$$
are bases for $\lin(P(G_1))$ and $\lin(P(G_2))$, respectively. Observe that $|B_i|=1+\dim(P(G_i))$ for $i\in [2]$. These notations and objects are fixed throughout this section. We proceed to state several preliminary claims, and present proofs for the non-trivial, non-routine statements.

\begin{custom}\label{composition}
{\bf Composition along separating cuts.} 
Given vectors $x\in \cR^{E(G_1)}$ and $y\in \cR^{E(G_2)}$ such that $x_e=y_e$ for all $e\in C$, we define \b{the composition} $z:=x\odot y\in \cR^E$ as follows: $z_e:=x_e$ if $e\in E(G_1)\setminus E(G_2)$, $z_e:=y_e$ if $e\in E(G_2)\setminus E(G_1)$, and $z_e:=x_e=y_e$ if $e\in C$. 
For each $e\in C$, let $I_e:=\{i\in [\hat{d}_1] : x^i_e=1\}$ and $J_e:=\{j\in [\hat{d}_2] : y^j_e=1\}$. Then both $I_e,J_e$ are nonempty as $G_1,G_2$ are matching-covered. Write $I_e=\{i_1,\ldots,i_k\}$ and $J_e=\{j_1,\ldots,j_\ell\}$, and let $$
			z^e_t:=x^{i_1} \odot y^{j_t},\quad t=1,\ldots,\ell
			\quad\text{and}\quad
			z^e_{\ell+t}:=x^{i_{1+t}}\odot y^{j_1},\quad t=1,\ldots,k-1.
		$$
	\rb{Let}{We define the \emph{composition}} $B_1\odot B_2:=\{z^e_i:e\in C, 1\leq i\leq |I_e|+|J_e|-1\}$. This is also known as the \emph{merger operation} (\cite{LM24}, \S6.3.1). The following statements can be readily checked for $B:=B_1\odot B_2$\rt{ (see the appendix 
	for a proof)}. \begin{enumerate}
\item[(1)] $B\subseteq P(G;C)\cap \{0,1\}^{E}$, $|B| = |B_1|+|B_2|-|C|$, and $B$ is a basis for $\lin(P(G;C))$. 
\item[(2)] If $B_i$ is a lattice basis for $L(G_i)$ for $i\in [2]$, then $B$ is a lattice basis for the lattice generated by the integral points in $P(G;C)$. 
\item[(3)] If $B_i$ is an integral basis for $\lin(P(G_i))$ for $i\in [2]$, then $B$ is an integral basis for $\lin(P(G;C))$.
\end{enumerate} Furthermore, we have the following. \begin{enumerate}
\item[(4)] Suppose $B_1\setminus \{x^1\}\subseteq \{x:x(D)=1\}$ for some $D\subseteq E(G_1)$, and $B_1\setminus \{x^1\}\not\subseteq \{x:x_f=0\}$ for any $f\in E(G_1)$. We can then apply the composition procedure in such a way that for some $z^\star\in B$, we have $B\setminus \{z^\star\}\subseteq \{z:z(D)=1\}$, $B\setminus \{z^\star\}\not\subseteq \{z:z_f=0\}$ for any $f\in E$, and $z^\star(D) = x^1(D)$. 
	 
	 To achieve this, for the element $e\in C$ such that $x^1_e=1$, we index $I_e$ such that $x^{i_1}\neq x^1$, which is possible as $B_1\setminus \{x^1\}\not\subseteq \{x:x_e=0\}$ and so $k\geq 2$. This guarantees that $x^1$ is composed only once with some other vector $y^j$ in $B_2$. We set $z^\star:= x^{1}\odot y^j$. \b{Note that $y^j$ must be composed with at least one other vector in $B_1\setminus \{x^1\}$, namely $x^{i_1}$, thereby guaranteeing that $B\setminus \{z^\star\}\not\subseteq \{z:z_f=0\}$ for any $f\in E$.}
	 \end{enumerate}
\end{custom}
		
\begin{custom}\label{uncrossing}
{\bf Uncrossing odd cuts.} Let $C_1=\delta(X_1),C_2=\delta(X_2)$ be separating cuts where $|X_1\cap X_2|$ is odd, and $X_1,X_2$ \emph{cross}, meaning that $X_1\cap X_2,X_1\setminus X_2,X_2\setminus X_1,\overline{X_1\cup X_2}$ are nonempty. Suppose further that $P(G;C_1)\cap P(G;C_2)\not\subseteq \{x:x_e=0\}$ for any edge $e\in E$. In particular, no edge of $G$ goes from $X_1\setminus X_2$ to $X_2\setminus X_1$, as such an edge would not belong to any perfect matching intersecting both $C_1,C_2$ exactly once. Let $I:=\delta(X_1\cap X_2)$ and $U:=\delta(X_1\cup X_2)$, both of which are odd cuts of $G$. As there is no edge from $X_1\setminus X_2$ to $X_2\setminus X_1$, it follows that $x(C_1) + x(C_2) = x(I) + x(U)$. Thus, $P(G;C_1)\cap P(G;C_2) = P(G;I)\cap P(G;U)$. 
\end{custom}

\begin{custom}\label{BvN-cut-contraction}
{\bf BvN cut-contractions.} 
\b{({\bf 1})} Suppose both $G_1,G_2$ are BvN. Then every facet of the polytope $P(G;C)$ is exposed by a non-negativity inequality. To see this, note that by definition, $P(G_1),P(G_2)$ are described by non-negativity inequalities and degree equations. Thus, it can be readily checked that $P(G;C)$ is described by non-negativity inequalities, degree equations, and $x(C)=1$, in turn implying the claim.
\b{{\bf (2)} If $G$ is a near-brick, and $C$ is a tight cut, then one of $G_1,G_2$, say $G_1$, is bipartite. In this case, $x(C)-1$ is a linear combination of the degree equations in $G\setminus \bar{X}$, so by (1), $G_1,G_2$ are BvN if and only if $G$ is BvN.}
\end{custom}

\begin{custom}\label{brick-count}
{\bf Number of bricks in cut-contractions.} Suppose $x(C)\geq 1$ exposes a face of $P(G)$ of dimension $d-i$, for some integer $i\geq 0$. Then $1+d-i=|B| =|B_1|+|B_2|-|C|$ \b{by \Cref{composition}, part (1)}, so $$
	|E|-|V|+2-b(G)-i = |E(G_1)|-|V(G_1)| +2 -b(G_1) + |E(G_2)|-|V(G_2)| +2 -b(G_2)-|C|,
	$$ implying in turn that $b(G_1)+b(G_2) = b(G)+ i$.
\end{custom}

\begin{custom}\label{brick-count-near-brick}
{\bf Number of bricks in cut-contractions of near-bricks.} {\bf (1)} In general, if one of $G_1,G_2$ is bipartite, then a simple counting argument implies that $C$ is a tight cut. That is, if $C$ is not tight, then $\min\{b(G_1),b(G_2)\}\geq 1$. {\bf (2)} The converse also holds if $G$ is a near-brick, by \Cref{brick-count}. {\bf (3)} Subsequently, if $G$ is a near-brick, then $C$ is facet-defining if, and only if, both $G_1,G_2$ are near-bricks.
\end{custom}

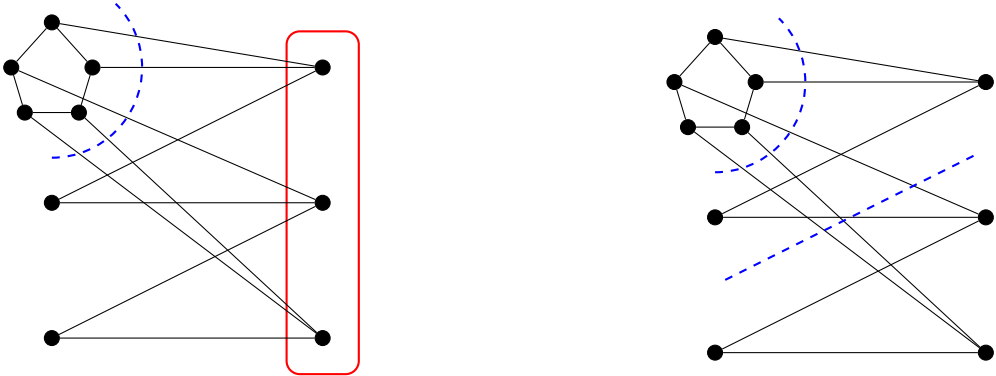
\begin{figure}[htbp]
    \centering
    \begin{minipage}{0.45\textwidth}
        \centering
        \begin{tikzpicture}[scale=1.2, mynode/.style={circle, draw, fill=black, inner sep=2pt}, cut/.style={blue, thick, dashed}]

\node[mynode] (a) at (0,4) {};  
\node[mynode] (b) at (0.45,3.5) {};
\node[mynode] (c) at (0.3,3) {};
\node[mynode] (d) at (-0.3,3) {};
\node[mynode] (e) at (-0.45,3.5) {};
\node[mynode] (f) at (0,2) {};
\node[mynode] (g) at (0,0.5) {}; 

\node[mynode] (x) at (3,3.5) {};
\node[mynode] (y) at (3,2) {};
\node[mynode] (z) at (3,0.5) {};

\begin{scope}
\draw[red, thick, rounded corners=5pt]
($(z) + (-0.4,-0.4)$) rectangle ($(x) + (0.4,0.4)$);
\end{scope}

\draw (a) -- (b);
\draw (b) -- (c);
\draw (c) -- (d);
\draw (d) -- (e);
\draw (e) -- (a);
\draw (a) -- (x);
\draw (b) -- (x);
\draw (c) -- (z);
\draw (d) -- (z);
\draw (e) -- (y);
\draw (f) -- (x);
\draw (f) -- (y);
\draw (g) -- (y);
\draw (g) -- (z);

\draw[cut] 
(0,2.5) arc (-90:45:1);

\end{tikzpicture}
    \end{minipage}
    \hfill
    \begin{minipage}{0.45\textwidth}
        \centering
        \begin{tikzpicture}[scale=1.2, mynode/.style={circle, draw, fill=black, inner sep=2pt}, cut/.style={blue, thick, dashed}]

\node[mynode] (a) at (0,4) {};  
\node[mynode] (b) at (0.45,3.5) {};
\node[mynode] (c) at (0.3,3) {};
\node[mynode] (d) at (-0.3,3) {};
\node[mynode] (e) at (-0.45,3.5) {};
\node[mynode] (f) at (0,2) {};
\node[mynode] (g) at (0,0.5) {}; 

\node[mynode] (x) at (3,3.5) {};
\node[mynode] (y) at (3,2) {};
\node[mynode] (z) at (3,0.5) {};

\draw (a) -- (b);
\draw (b) -- (c);
\draw (c) -- (d);
\draw (d) -- (e);
\draw (e) -- (a);
\draw (a) -- (x);
\draw (b) -- (x);
\draw (c) -- (z);
\draw (d) -- (z);
\draw (e) -- (y);
\draw (f) -- (x);
\draw (f) -- (y);
\draw (g) -- (y);
\draw (g) -- (z);

\node (P) at (0,1.25) {};
\node (Q) at (3,2.75) {};

\draw[cut] (P) -- (Q);
\draw[cut] (0,2.5) arc (-90:45:1);

\end{tikzpicture}
    \end{minipage}
    \caption{Barrier and equivalent cuts. (Left): an example of a near-brick with a tight cut (blue dashed). The corresponding barrier is circled in red. (Right): The two dashed blue cuts are equivalent separating cuts in a near-brick. Contracting two disjoint shores results in $4$-cycle, a bipartite graph.}
    \label{fig:tikz_side_simple}
\end{figure}

\begin{custom}\label{tight-cut-barrier}
{\bf Tight cuts in near-bricks and barriers.} Suppose $G$ is a near-brick, and $C$ is a tight cut of $G$. Then one of $G_1,G_2$ is a near-brick while the other one, say $G_1$, is bipartite, by \Cref{brick-count}. Let $R\subseteq V$ be the part in a bipartition of $G_1$ where the shrunk vertex $x\notin R$. Observe that $R$ is an independent set of $G$, and $G\setminus R$ has exactly $|R|$ connected components, one of which is $X$, and all others are singletons. We say that $R$ is a \emph{barrier} in $G$, and $C$ is a \emph{barrier cut} with barrier $R$. See Figure~\ref{fig:tikz_side_simple} (left) for an illustrated example. A barrier is \emph{maximal} if it is not contained in another barrier.\footnote{Barriers are defined more broadly and for all matching-covered graphs, but for our purposes, we shall focus only on the (special type of) barriers defined here for near-bricks.}
\end{custom}

\begin{custom}\label{tight-cut-fdi}
{\bf Tight cuts in near-bricks and facet-defining cuts.} Suppose $G$ is a near-brick, $C$ is a tight cut of $G$, $G_1$ is bipartite, and $G_2$ is a non-BvN near-brick. If $D=\delta(Y)$ is a separating facet-defining cut of $G_2$ with $Y\subset X$, then $D$ is a separating facet-defining cut of $G$. 

To this end, as $D$ is separating in $G_2$, and $C$ is separating in $G$, it follows that $D$ is separating in~$G$ \b{as we remarked already in the introduction after defining separating cuts}. 
As $G_2$ is a near-brick, both $G_2/Y,G_2/\bar{Y}$ are near-bricks, by \Cref{brick-count-near-brick}, part (3). Thus, $G/\bar{Y}=G_2/\bar{Y}$ is a near-brick. As $G_2/Y$ 
is a near-brick, $C=\delta(\bar{X})$ is a tight cut in $G/Y$, and $G_1$ is bipartite, it follows that $G/Y$ is a near-brick. Subsequently, both $G/Y,G/\bar{Y}$ are near-bricks, implying that $D$ is a separating facet-defining cut of $G$, by \Cref{brick-count-near-brick}, part (3).
\end{custom}

\begin{custom}\label{equivalent}
{\bf Equivalent cuts in near-bricks.} See Figure~\ref{fig:tikz_side_simple} (right) for an illustrated example helpful for the following arguments. Two odd cuts $C_1,C_2$ are \emph{equivalent} if $x(C_1)=x(C_2)$ for all $x\in P(G)$. Suppose $G$ is a near-brick, and $C_1=\delta(X_1),C_2=\delta(X_2)$ are separating cuts that define one and the same facet, where $|X_1\cap X_2|$ is odd \b{which can be assumed without loss of generality by replacing one of the sets with its complement}. We claim that $C_1,C_2$ are equivalent; furthermore, if $X_1\subset X_2$, then $G/X_1/\bar{X}_2$ is bipartite. To this end, note that $G/\bar{X}_1,G/\bar{X}_2$ are near-bricks, \b{by \ref{brick-count-near-brick}, part (3)}.

If $X_1=X_2$, then the claim is clear.

If $X_1,X_2$ do not cross, say $X_1\subset X_2$, then $\delta(X_1)$ is a tight cut in the near-brick $G/\bar{X}_2$, and as $G/\bar{X}_1$ is a near-brick, it follows that $G/X_1/\bar{X}_2$ is bipartite matching-covered, with $x_1,\bar{x}_2$ on opposite sides of any bipartition \b{(as $G$ has a perfect matching $M$ such that $|M\cap \delta(X_1)|,|M\cap \delta(X_2)|>1$)}. A simple counting argument now implies that for every perfect matching $M$ of $G$, $|M\cap \delta(X_1)|=|M\cap \delta(X_2)|$, so $C_1,C_2$ are equivalent cuts in $G$.

Otherwise, $X_1,X_2$ cross. Let $I:=\delta(X_1\cap X_2)$ and $U:=\delta(X_1\cup X_2)$. Given that $P(G;C_1)\cap P(G;C_2) = P(G;C_1)\not\subseteq \{x:x_e=0\}$ for any edge $e\in E$, it follows from \Cref{uncrossing} that $P(G;C_1) = P(G;I)\cap P(G;U)$. As $P(G;C_1)$ is a facet of $P(G)$, at least one of $I,U$ must define the same facet as $C_1$, say it is $I$. By applying the argument above to $X_1,X_1\cap X_2$, and also to $X_2,X_1\cap X_2$, we obtain that $C_1,I$ and $C_2,I$ are equivalent, so $C_1,C_2$ are equivalent, as claimed.
\end{custom}

\begin{custom}\label{linear-basis-pete}
{\bf Basis for near-bricks with a Petersen brick.} 
Suppose $G$ is a near-brick with a Petersen brick $H$ obtained through a tight cut decomposition. Let $Y$ be the vertex set of a $5$-cycle of $H$, and let $D:=\delta(Y)$. Then there exist perfect matchings $M_0,M_1,\ldots,M_d$ of $G$ whose incidence vectors form a basis for $\lin(P(G))$, where $|M_0\cap D|=5$, $|M_i\cap D|=1$ for all $i\in [d]$, and every edge of $G$ is contained in at least one of $M_1,\ldots,M_d$.

To see this, note first that the Petersen graph $\pete$ has exactly six perfect matchings $N_0,N_1,\ldots,N_5$, whose incidence vectors form a basis for $\lin(P(\pete))$. One of the perfect matchings, say $N_0$, intersects $\delta_{\pete}(Y)$ in five edges, while the remaining perfect matchings intersect $\delta_{\pete}(Y)$ just once. Note further that every edge of $\pete$ belongs to exactly two perfect matchings, so at least one of these must be among $N_1,\ldots,N_5$.

Subsequently, given that $H$ is obtained from $\pete$ by adding some $p$ edges parallel to existing edges, we can find perfect matchings $N_6,\ldots,N_{p+5}$ of $H$ such that the incidence vectors of $N_0,N_1,\ldots,N_{p+5}$ form a basis for $\lin(P(H))$, $|N_0\cap D|=5$, and $|N_i\cap D|=1$ for all $i\in [p]$, and every edge of $H$ belongs to at least one of $N_1,\ldots,N_{p+5}$.

The existence of $M_0,M_1,\ldots,M_d$ now follows by recursively applying \Cref{composition}, part (4).
\end{custom}

\begin{custom}\label{cut-face-triples}
{\bf Cut and face triples.}
A \emph{cut triple} is a tuple $(C_1,C_2,C_3)$ of odd cuts of $G$ of the form $C_i=\delta(X_i), i\in [3]$ such that $X_1\subset X_2\subset X_3$. A \emph{face triple} is a tuple $(F_1,F_2,F_3)$ of faces of $P(G)$ such that $F_2\cap F_1 = F_2\cap F_3$, $F_2\cap F_3$ is not contained in $\{x:x_e=0\}$ for any $e\in E$, and at least one of $F_2\setminus F_1,F_2\setminus F_3$ is nonempty. We claim that a cut triple cannot define a face triple, that is, we cannot have $P(G;C_i)=F_i,i=1,2,3$.

Suppose otherwise. By symmetry between $C_1$ and $C_3$, we may assume that $F_2\setminus F_1\neq \emptyset$. Let $M$ be a perfect matching such that $|M\cap C_2|=1$ and $|M\cap C_1|>1$. Write $M\cap C_2=\{f\}$. As $F_2\cap F_3\not\subseteq \{x:x_f=0\}$, there exists a perfect matching $M'$ such that $M'\cap C_2=\{f\}$ and $|M'\cap C_3|=1$.

Let $M''$ be the perfect matching such that $M''\cap C_2=\{f\}$, and agrees with $M$ in $G[X_2]$ and with $M'$ in $G[\bar{X_2}]$. Then $|M''\cap C_1|=|M\cap C_1|>1$ and $|M''\cap C_3|=|M'\cap C_3|=1$. Subsequently, the incidence vector of $M''$ belongs to $F_2\cap F_3$ but not $F_2\cap F_1$, a contradiction as $F_2\cap F_1= F_2\cap F_3$.
\end{custom}

\section{Proof of the Intersection Theorem}\label{sec:intersection-theorem}

Let $G=(V,E)$ be a Petersen-free non-BvN near-brick. We prove by induction on $|E|$ that there exist a perfect matching $M$ and a separating facet-defining cut $C$ such that $|M\cap C|=3$; we shall call $(M,C)$ a \emph{desired pair}. We proceed in four stages.

\paragraph{From near-bricks to bricks.} We first reduce to the case where $G$ is a brick. To this end, suppose $G$ has a tight cut $\delta(X)$, and let $G_1,G_2$ be the corresponding cut-contractions, where $G_1$ is bipartite, and $G_2$ is a Petersen-free near-brick. We know that $G_2$ is non-BvN, by \Cref{BvN-cut-contraction}. By the induction hypothesis, $G_2$ has a perfect matching $N$ and a separating facet-defining cut $C$ such that $|N\cap C|=3$. Then $C$ is a separating facet-defining cut of $G$, by \Cref{tight-cut-fdi}. Thus, by extending $N$ to a perfect matching $M$ of $G$, we see that $(M,C)$ is a desired pair. 

\paragraph{The facets.} Next we reduce to the case where no facet of $P(G)$ is exposed by a non-negativity inequality. 

Suppose some facet of $P(G)$ is exposed by a non-negativity inequality. As $G$ is \b{is a non-BvN brick}, there exists one such facet, say exposed by $x_e\geq 0$, that is adjacent to a facet not exposed by a non-negativity inequality. 

We claim that $G\setminus e$ is a non-BvN near-brick. To this end, let $F$ be the union of all perfect matchings of $G\setminus e$. Then \b{the restriction graph} $G|F \b{:=(V,F)}$ is a matching-covered graph such that $P(G|F)$ is obtained from the facet $P_e(G)$ after removing the coordinates that are equal to $0$ throughout the facet. Our choice of the facet $P_e(G)$ implies that $P(G|F)$ has a facet not exposed by a non-negativity inequality, so $G|F$ is not BvN, implying in turn that $G|F$ is not bipartite, so $b(G|F)\geq 1$. On the other hand, 
$$d-1=\dim(P(G|F)) = |F|-|V|+1-b(G|F) = d +1 -|E\setminus F| - b(G|F),$$ implying in turn that $|E\setminus F|+b(G|F) = 2$, so $|E\setminus F|=b(G|F)=1$. Thus, $G|F=G\setminus e$ is a non-BvN near-brick. 

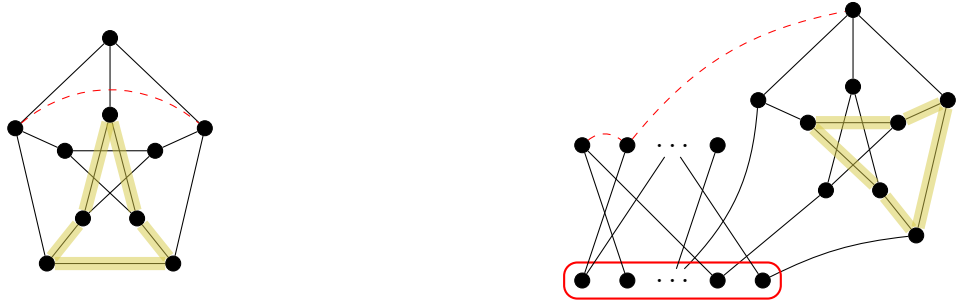
\begin{figure}[htbp]
    \centering
    \begin{minipage}{0.45\textwidth}
        \centering
        \begin{tikzpicture}[scale=1.2, mynode/.style={circle, draw, fill=black, inner sep=2pt}, cut/.style={blue, thick, dashed}, spoke/.style={line width=5pt, color=yellow!80!black, opacity=0.5}]

\node[mynode] (p1) at (0,2.5) {};
\node[mynode] (p2) at (-1.05,1.5) {};
\node[mynode] (p3) at (-0.7,0) {};
\node[mynode] (p4) at (0.7,0) {};
\node[mynode] (p5) at (1.05,1.5) {};
\node[mynode] (p6) at (0,1.65) {};
\node[mynode] (p7) at (-0.5,1.25) {};
\node[mynode] (p8) at (-0.3,0.5) {};
\node[mynode] (p9) at (0.3,0.5) {};
\node[mynode] (p10) at (0.5,1.25) {};

\draw (p1)--(p2);
\draw (p3)--(p2);
\draw (p4)--(p3);
\draw (p5)--(p4);
\draw (p1)--(p5);
\draw (p1)--(p6);
\draw (p2)--(p7);
\draw (p3)--(p8);
\draw (p4)--(p9);
\draw (p5)--(p10);
\draw (p6)--(p8);
\draw (p8)--(p10);
\draw (p7)--(p9);
\draw (p9)--(p6);
\draw (p10)--(p7);

\draw[color=red, dashed] (p2) to[bend left = 40] (p5);

\draw[spoke] (p6)--(p8);
\draw[spoke] (p8)--(p3);
\draw[spoke] (p3)--(p4);
\draw[spoke] (p4)--(p9);
\draw[spoke] (p9)--(p6);

\end{tikzpicture}
    \end{minipage}
    \hfill
    \begin{minipage}{0.45\textwidth}
        \centering
        \begin{tikzpicture}[scale=1.2, mynode/.style={circle, draw, fill=black, inner sep=2pt}, cut/.style={blue, thick, dashed}, spoke/.style={line width=5pt, color=yellow!80!black, opacity=0.5}]

\node[mynode] (p1) at (0,2.5) {};
\node[mynode] (p2) at (-1.05,1.5) {};
\node[mynode] (p4) at (0.7,0) {};
\node[mynode] (p5) at (1.05,1.5) {};
\node[mynode] (p6) at (0,1.65) {};
\node[mynode] (p7) at (-0.5,1.25) {};
\node[mynode] (p8) at (-0.3,0.5) {};
\node[mynode] (p9) at (0.3,0.5) {};
\node[mynode] (p10) at (0.5,1.25) {};

\node[mynode] (a1) at (-3,1) {};
\node[mynode] (a2) at (-2.5,1) {};
\node         (a3) at (-2,1) {$\ldots$};
\node[mynode] (a4) at (-1.5,1) {};
\node[mynode] (b1) at (-3,-0.5) {};
\node[mynode] (b2) at (-2.5,-0.5) {};
\node         (b3) at (-2,-0.5) {$\ldots$};
\node[mynode] (b4) at (-1.5,-0.5) {};
\node[mynode] (b5) at (-1,-0.5) {};

\draw (p1)--(p2);
\draw (p2) to[bend left = 20] (b3);
\draw (p4) to[bend right = 10] (b5);
\draw (p5)--(p4);
\draw (p1)--(p5);
\draw (p1)--(p6);
\draw (p2)--(p7);
\draw (p8) to[bend left = 0] (b4);
\draw (p4)--(p9);
\draw (p5)--(p10);
\draw (p6)--(p8);
\draw (p8)--(p10);
\draw (p7)--(p9);
\draw (p9)--(p6);
\draw (p10)--(p7);

\begin{scope}
  \draw[red, thick, rounded corners=5pt]
    ($(b1) + (-0.2,-0.2)$) rectangle ($(b5) + (0.2,0.2)$);
\end{scope}

\draw (a1)--(b2);
\draw (a1)--(b4);
\draw (a2)--(b1);
\draw (a3)--(b5);
\draw (a4)--(b3);
\draw (a3)--(b1);

\draw[color=red, dashed] (a1) to[bend left = 40] (a2);
\draw[color=red, dashed] (a2) to[bend left = 20] (p1);
\draw[spoke] (p7)--(p9);
\draw[spoke] (p9)--(p4);
\draw[spoke] (p5)--(p4);
\draw[spoke] (p5)--(p10);
\draw[spoke] (p10)--(p7);

\end{tikzpicture}
    \end{minipage}
    \caption{Proof of the Intersection Theorem. (Left): Case $1$. the edge $e$ (red dashed) connects two vertices at distance two. The $5$-cycle of $H$ that we could pick is highlighted. (Right): Case $2$. The only maximal barrier is shown by the red box and the options for $e$ are shown as red dashed edges. The $5$-cycle of $H$ that we could pick is highlighted; its choice depends on the perfect matching $M$ (not shown).}
    \label{fig:tikz_side_simple2}
\end{figure}

\begin{figure}[htbp]
    \centering
    \begin{minipage}{0.45\textwidth}
        \centering
        \begin{tikzpicture}[scale=1.0, mynode/.style={circle, draw, fill=black, inner sep=2pt}, cut/.style={blue, thick, dashed}, spoke/.style={line width=5pt, color=yellow!80!black, opacity=0.5}]

\node[mynode] (p1) at (0,2.5) {};
\node[mynode] (p2) at (-1.05,1.5) {};
\node[mynode] (p5) at (1.05,1.5) {};
\node[mynode] (p6) at (0,1.65) {};
\node[mynode] (p7) at (-0.5,1.25) {};
\node[mynode] (p8) at (-0.3,0.5) {};
\node[mynode] (p9) at (0.3,0.5) {};
\node[mynode] (p10) at (0.5,1.25) {};

\node[mynode] (a1) at (-3,1) {};
\node[mynode] (a2) at (-2.5,1) {};
\node         (a3) at (-2,1) {$\ldots$};
\node[mynode] (a4) at (-1.5,1) {};
\node[mynode] (b1) at (-3,-0.5) {};
\node[mynode] (b2) at (-2.5,-0.5) {};
\node         (b3) at (-2,-0.5) {$\ldots$};
\node[mynode] (b4) at (-1.5,-0.5) {};
\node[mynode] (b5) at (-1,-0.5) {};
\node[mynode] (c4) at (3,1) {};
\node[mynode] (c3) at (2.5,1) {};
\node         (c2) at (2,1) {$\ldots$};
\node[mynode] (c1) at (1.5,1) {};
\node[mynode] (d5) at (3,-0.5) {};
\node[mynode] (d4) at (2.5,-0.5) {};
\node         (d3) at (2,-0.5) {$\ldots$};
\node[mynode] (d2) at (1.5,-0.5) {};
\node[mynode] (d1) at (1,-0.5) {};

\draw (p1)--(p2);
\draw (p2) to[bend left = 20] (b3);
\draw (d1) to[bend right = 20] (b5);
\draw (p5) to[bend right = 20] (d3);
\draw (p1)--(p5);
\draw (p1)--(p6);
\draw (p2)--(p7);
\draw (p8) to[bend left = 0] (b4);
\draw (p9) to[bend right = 0] (d2);
\draw (p5)--(p10);
\draw (p6)--(p8);
\draw (p8)--(p10);
\draw (p7)--(p9);
\draw (p9)--(p6);
\draw (p10)--(p7);

\begin{scope}
  \draw[red, thick, rounded corners=5pt]
    ($(b1) + (-0.2,-0.2)$) rectangle ($(b5) + (0.2,0.2)$);
\end{scope}
\begin{scope}
  \draw[red, thick, rounded corners=5pt]
    ($(d1) + (-0.2,-0.2)$) rectangle ($(d5) + (0.2,0.2)$);
\end{scope}

\draw (a1)--(b2);
\draw (a1)--(b4);
\draw (a2)--(b1);
\draw (a3)--(b5);
\draw (a4)--(b3);
\draw (a3)--(b1);

\draw (d1)--(c2);
\draw (c1)--(d4);
\draw (d2)--(c1);
\draw (c3)--(d5);
\draw (c4)--(d3);
\draw (d3)--(c1);

\draw[color=red, dashed] (a1) to[bend left = 40] (c4);
\draw[spoke] (p7)--(p9);
\draw[spoke] (p9)--(p6);
\draw[spoke] (p6)--(p8);
\draw[spoke] (p8)--(p10);
\draw[spoke] (p10)--(p7);

\end{tikzpicture}
    \end{minipage}
    \hfill
    \begin{minipage}{0.45\textwidth}
        \centering
        \begin{tikzpicture}[scale=1.0, mynode/.style={circle, draw, fill=black, inner sep=2pt}, cut/.style={blue, thick, dashed}, spoke/.style={line width=5pt, color=yellow!80!black, opacity=0.5}]

\node[mynode] (p1) at (0,2.5) {};
\node[mynode] (p2) at (-1.05,1.5) {};
\node[mynode] (p4) at (0.7,0) {};
\node[mynode] (p6) at (0,1.65) {};
\node[mynode] (p7) at (-0.5,1.25) {};
\node[mynode] (p8) at (-0.3,0.5) {};
\node[mynode] (p9) at (0.3,0.5) {};
\node[mynode] (p10) at (0.5,1.25) {};

\node[mynode] (a1) at (-3,1) {};
\node[mynode] (a2) at (-2.5,1) {};
\node         (a3) at (-2,1) {$\ldots$};
\node[mynode] (a4) at (-1.5,1) {};
\node[mynode] (b1) at (-3,-0.5) {};
\node[mynode] (b2) at (-2.5,-0.5) {};
\node         (b3) at (-2,-0.5) {$\ldots$};
\node[mynode] (b4) at (-1.5,-0.5) {};
\node[mynode] (b5) at (-1,-0.5) {};
\node[mynode] (c4) at (3,1) {};
\node[mynode] (c3) at (2.5,1) {};
\node         (c2) at (2,1) {$\ldots$};
\node[mynode] (c1) at (1.5,1) {};
\node[mynode] (d5) at (3,-0.5) {};
\node[mynode] (d4) at (2.5,-0.5) {};
\node         (d3) at (2,-0.5) {$\ldots$};
\node[mynode] (d2) at (1.5,-0.5) {};
\node[mynode] (d1) at (1,-0.5) {};

\draw (p1)--(p2);
\draw (p2) to[bend left = 20] (b3);
\draw (p4) to[bend right = 20] (b5);
\draw (p4) to[bend left = 10] (d1);
\draw (p1) to[bend right = 0] (d3);
\draw (p1)--(p6);
\draw (p2)--(p7);
\draw (p8) to[bend left = 0] (b4);
\draw (p9) to[bend right = 0] (p4);
\draw (p10) to[bend right = 0] (d2);
\draw (p6)--(p8);
\draw (p8)--(p10);
\draw (p7)--(p9);
\draw (p9)--(p6);
\draw (p10)--(p7);

\begin{scope}
  \draw[red, thick, rounded corners=5pt]
    ($(b1) + (-0.2,-0.2)$) rectangle ($(b5) + (0.2,0.2)$);
\end{scope}
\begin{scope}
  \draw[red, thick, rounded corners=5pt]
    ($(d1) + (-0.2,-0.2)$) rectangle ($(d5) + (0.2,0.2)$);
\end{scope}

\draw (a1)--(b2);
\draw (a1)--(b4);
\draw (a2)--(b1);
\draw (a3)--(b5);
\draw (a4)--(b3);
\draw (a3)--(b1);

\draw (d1)--(c2);
\draw (c1)--(d4);
\draw (d2)--(c1);
\draw (c3)--(d5);
\draw (c4)--(d3);
\draw (d3)--(c1);

\draw[color=red, dashed] (a1) to[bend left = 40] (c4);
\draw[spoke] (p7)--(p9);
\draw[spoke] (p9)--(p6);
\draw[spoke] (p6)--(p8);
\draw[spoke] (p8)--(p10);
\draw[spoke] (p10)--(p7);

\end{tikzpicture}
    \end{minipage}
    \caption{Proof of the Intersection Theorem: Two options for Case $3$. Two maximal barriers are in red and $e$ is shown in dashed. A $5$-cycle of $H$ that we could pick is highlighted; its choice depends on the perfect matching $M$ (not shown). (Left): the contraction vertices of two maximal barriers are adjacent in $H$. (Right): the two contraction vertices are at distance two.}
    \label{fig:tikz_side_simple_case3}
\end{figure}
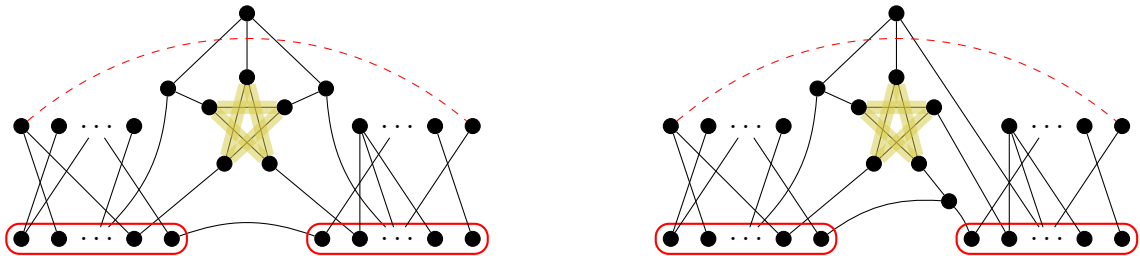

\medskip
Suppose $G\setminus e$ has a Petersen brick. Let $H$ be the Petersen brick of $G\setminus e$, obtained through a tight cut decomposition of $G\setminus e$. We have the following three cases, which come with illustrations in Figures~\ref{fig:tikz_side_simple2} and~\ref{fig:tikz_side_simple_case3}. \begin{enumerate}
\item[Case 1.] $H=G\setminus e$. In this case, as $G$ is a non-Petersen brick, $e$ must connect a pair of vertices of $H$ at distance $2$. 
\end{enumerate} Otherwise, $H\neq G\setminus e$, so $G\setminus e$ has a tight cut. Every tight cut of $G\setminus e$ is a barrier cut of the form $\delta_{G\setminus e}(X)$ with barrier $B\subset \bar{X}$, by \Cref{tight-cut-barrier}. Observe that $P_e(G)\subseteq P(G;\delta_G(X))$, and since $P(G;\delta_G(X))\neq P(G)$ as $G$ is a brick, \rb{we must have}{dimension counting implies} $P_e(G)= P(G;\delta_G(X))$. A simple counting argument implies that $e$ joins either two singleton components of $G\setminus e\setminus B$, or it joins a singleton component to $X$. This implies that there exist either one or two maximal barriers only, summarized as follows.
\begin{enumerate}
\item[Case 2.] $H=G\setminus e/\bar{X}$, where $\delta_{G\setminus e}(X)$ is a barrier cut with barrier $B\subset \bar{X}$, and $B$ is the unique maximal barrier of $G\setminus e$. 
\item[Case 3.] $H=G\setminus e/\bar{X}_1/\bar{X}_2$, where $\delta_{G\setminus e}(X_i)$ is a barrier cut with barrier $B_i\subset \bar{X}_i$ for $i\in [2]$, $\bar{X}_1\cap \bar{X}_2=\emptyset$, and $e$ joins a singleton component of $G\setminus e\setminus B_1$ to a singleton component of $G\setminus e\setminus B_2$. Furthermore, $B_1,B_2$ are the only maximal barriers of $G\setminus e$. 
\end{enumerate}

It can be shown through an elementary though terse argument that in each of Cases 1-3, there exist a $5$-cycle of $H$ with vertex set $Y$, and perfect matchings $M,M'$ of $G$ containing $e$, such that $Y$ does not use contraction vertices of $H$, $Y$ is not incident to either ends of $e$, $|M\cap \delta_G(Y)|=3$, and $|M'\cap \delta_G(Y)|=1$.\footnote{For a detailed analysis, we refer the reader to the proof of (\cite{CLM1}, Theorem 5.4), or (\cite{LM24}, Theorem 13.6).} Moving on, we claim that $\delta_G(Y)$ is a separating facet-defining cut of $G$. Clearly, $\delta_{G\setminus e}(Y)$ is separating in $G\setminus e$ as it is so in $H$. Thus, as $e\in M'$ and $|M'\cap \delta_G(Y)|=1$, it follows that $\delta_G(Y)$ is separating in~$G$. To see that $\delta_G(Y)$ is facet-defining in $G$, pick perfect matchings $M_1,\ldots,M_{d}$ whose incidence vectors lie inside the facet $P_e(G)$ and are linearly independent; we may pick the first $d-1$ perfect matchings to have intersection $1$ with $\delta_G(Y)$, and the last one with intersection $5$, by \Cref{linear-basis-pete}. We then swap $M_d$ with the perfect matching $M'\ni e$ which satisfies $|M'\cap \delta(Y)|=1$, in order to obtain $d$ linearly independent vectors $M_1,\ldots,M_{d-1},M'$ in the face $P(G;\delta_G(Y))$, implying in turn that $\delta_G(Y)$ defines a facet of $P(G)$. Consequently, $(M,\delta_G(Y))$ is a desired pair.\medskip

We may therefore assume that $G\setminus e$ is a Petersen-free non-BvN near-brick. By the induction hypothesis, there exist a perfect matching $M$ and a separating facet-defining cut $\delta_{G\setminus e}(X)$ of $G\setminus e$ such that $|M\cap \delta_{G\setminus e}(X)|=3$. Let $C:=\delta_G(X)$. Observe that $P(G;C)$ is either facet of $P(G)$, or a $(d-2)$-face of $P(G)$ contained in $P_e(G)$. Observe further that $P(G;C)$ is not contained in $\{x:x_f=0\}$ for any $f\in E\setminus e$. 

Suppose first that $P(G;C)$ is a facet of $P(G)$. Given that $e\notin M$ and $|M\cap C|\neq 1$, it follows that $P(G;C)\neq P_e(G)$. Subsequently, $P(G;C)$ is not contained in $\{x:x_f=0\}$ for any edge $f\in E$, so $C$ is a separating facet-defining cut of $G$. Subsequently, $(M,C)$ is a desired pair. 

In the remaining case, $P(G;C)$ is a $(d-2)$-face of $P(G)$, contained in the facet $P_e(G)$. Consider now the other facet of $P(G)$ containing $P(G;C)$. It cannot be exposed by a non-negativity inequality, given that $P(G;C)$ is not contained in $\{x:x_f=0\}$ for any $f\in E\setminus e$. In particular, the other facet is exposed by a separating facet-defining cut of $G$, say $D$. In particular, $D\setminus e$ defines the same facet as $C\setminus e$ \b{in $P(G\setminus e)$}, and so $D\setminus e,C\setminus e$ must be equivalent cuts as $G\setminus e$ is a near-brick, by \Cref{equivalent}. Thus, $|M\cap D|=|M\cap C|=3$, implying in turn that $(M,D)$ is a desired pair. 

\paragraph{The $(d-2)$-faces.} Thus every non-negativity inequality defines a face of $P(G)$ of dimension at most $d-2$. In this stage, we reduce to the case where every $(d-2)$-face of $P(G)$ is exposed by a non-negativity inequality. 

Suppose some $(d-2)$-face of $P(G)$ is not exposed by a non-negativity inequality. Take a facet containing this face, which is exposed by a separating facet-defining cut $C=\delta(X)$. Then both $C$-contractions of $G$ are near-bricks, by \Cref{brick-count-near-brick}, part (3). Given that the polytope $P(G;C)$ has a facet not exposed by a non-negativity inequality, one of the $C$-contractions, say $G/\bar{X}$, must be non-BvN, by \Cref{BvN-cut-contraction}.

We may assume that $G/\bar{X}$ is a non-BvN brick. To argue this, note first that the unique brick of $G/\bar{X}$ is non-BvN, by \Cref{BvN-cut-contraction}. If $G/\bar{X}$ is not this brick, then choose a minimal vertex subset $Y\subset X$ such that $\delta(Y)$ is a tight cut of $G/\bar{X}$. In particular, $P(G;C)\subseteq P(G;\delta(Y))$. As $G$ is a brick, $\delta(Y)$ is not tight in $G$, so $P(G;C)= P(G;\delta(Y))$, therefore $C,\delta(Y)$ are equivalent cuts and $G/Y/\bar{X}$ is bipartite, by \Cref{equivalent}. This, together with our minimal choice of $Y$, implies that $G/\bar{X}/\bar{Y}=G/\bar{Y}$ is a brick. Subsequently, by changing $C=\delta(X)$ \b{to $\delta(Y)$} and without changing the corresponding facet, if necessary, we may assume $G/\bar{X}$ is a non-BvN brick.

If $G/\bar{X}$ is a Petersen brick, \b{the degree of $\bar{x}$ in the contraction is exactly three, so} any perfect matching $M$ intersecting $C$ more than once must intersect $C$ exactly three times, so $(M,C)$ is a desired pair. 

Therefore, we may assume that $G/\bar{X}$ is a non-Petersen non-BvN brick. By the induction hypothesis, $G/\bar{X}$ has a separating facet-defining cut $R$ and a perfect matching $N$ such that $|R\cap N|=3$. Extend $N$ to a perfect matching $M$ of $G$. Then $|R\cap M|=3$ and $|C\cap M|=1$. \b{Since $C$ is separating in $G$}, $R$ is \b{also} a separating cut of $G$. \rb{Observe that}{Since $C$ is facet-defining,} $P(G;R)$ is either a facet, or a $(d-2)$-face of $P(G)$ contained in $P(G;C)$.

If $P(G;R)$ is a facet of $P(G)$, then $(M,R)$ is a desired pair.

Otherwise, $P(G;R)$ is a $(d-2)$-face of $P(G)$ contained in $P(G;C)$. Let $P(G;D)$ be the other facet of $P(G)$ containing $P(G;R)$. Observe that $x(R)\geq 1$ and $x(D)\geq 1$ define the same facet of the polytope $P(G;C)$.  Observe further that $(P(G;R),P(G;C),P(G;R))$, as well as any permutation of $(P(G;R),P(G;C),P(G;D))$ without $P(G;R)$ in the middle, is a face triple. 

We shall prove that $(M,D)$ is a desired pair. To this end, write $D=\delta(Z)$ where $Z\subset V$ and $|X\cap Z|$ is odd. There are three cases:
 \begin{enumerate}
\item[Case 1.] Suppose $Z\subset X$. Then $D$ is an odd cut of $G/\bar{X}$, so $R,D$ define the same facet of $P(G/\bar{X})$, so $R,D$ are equivalent cuts of $G/\bar{X}$ as $G/\bar{X}$ is a brick, by \Cref{equivalent}. Thus, $|D\cap N|=|R\cap N|=3$, so $|D\cap M|=|D\cap N|=3$, therefore $(M,D)$ is a desired pair.
\item[Case 2.] Suppose $X\subset Z$. Then $(R,C,D)$ would be a cut triple defining a face triple, thus contradicting \Cref{cut-face-triples}.
\item[Case 3.] In the remaining case, $X$ and $Z$ cross. Let $I:=\delta(X\cap Z)$ and $U:=\delta(X\cup Z)$. 
Given that $P(G;C)\cap P(G;D) = P(G;R)\not\subseteq \{x:x_e=0\}$ for any $e\in E$, it follows from \Cref{uncrossing} that $x(C)+x(D) = x(I)+x(U)$ and $P(G;R)=P(G;C)\cap P(G;D) = P(G;I)\cap P(G;U)$. 

If $\{P(G;I),P(G;U)\}=\{P(G;R),P(G;C)\}$ or $\{P(G;R),P(G;D)\}$, then 
$(I,D,U)$ or $(I,C,U)$ would be a cut triple defining a face triple, respectively, thus contradicting \Cref{cut-face-triples}.

If $\{P(G;I),P(G;U)\}=\{P(G;R),P(G)\}$, then we have $P(G;I)=P(G;R)$ and $P(G;U)=P(G)$. This holds since otherwise $(R,C,U)$ would be a cut triple defining a face triple, thus contradicting \Cref{cut-face-triples}. The equality $P(G;U)=P(G)$ implies that $U$ is a trivial cut as $G$ is a brick, while $P(G;I)=P(G;R)$ implies that $R,I$ are equivalent cuts in the brick $G/\bar{X}$, so $x(I)=x(R)$ for all $x\in P(G)$. Consequently, $$x(D) = x(I)+x(U)-x(C)=x(R)+1-x(C),$$ so whenever $x(C)=1$, then $x(R)=x(D)$, so $|D\cap M|=|R\cap M|=3$, so $(M,D)$ is a desired pair. 

Otherwise, $\{P(G;I),P(G;U)\}=\{P(G;C),P(G;D)\}$, then we have $P(G;I)=P(G;D)$ and $P(G;U)=P(G;C)$. This holds since otherwise $(R,C,U)$ would be cut triple defining a face triple, thus contradicting \Cref{cut-face-triples}. Now by replacing $D$ with $I$ we fall back to Case 1.
\end{enumerate}

\paragraph{The finale.} We have reduced to the case where $G$ is a non-BvN non-Petersen brick where every $(d-2)$-face of $P(G)$ is exposed by a non-negativity inequality. The existence of a desired pair now follows from the Petersen Graph Lemma, thus finishing the proof.\qed

\section{Proof of the Main Theorem}\label{sec:integral-basis-theorem}

In this final section, we prove the Integral Basis Theorem, and then obtain the Main Theorem as a corollary. For a polyhedron $P\subseteq \cR^n$ and $k\geq 0$, define $kP$ as the set of all points of the form $\sum_{p\in P}\lambda_p p$ where $\lambda\in \cR^P_{\geq 0}$ and $\1^\top \lambda=k$. Here, $\1$ denotes the all-ones vector of appropriate dimension. $P$ has the \emph{integer decomposition property} if for every integer $k\geq 1$, every integral point in $kP$ can be written as the sum of $k$ integral points in $P$. We need the following preliminary.

\begin{theorem}[\cite{ACLS25+}]\label{IDP-integral-basis}
Let $P\subseteq \cR^n$ be a pointed polyhedron with the integer decomposition property, whose affine hull is of the form $\{x:Ax=b\}$ for $A\in \cZ^{m\times n},b\in \cZ^m$ such that $m \geq 1$, $b\neq \0$, and $\gcd\{b_i:i\in [m]\}=1$. Then $P\cap \cZ^n$ contains an integral basis for $\lin(P)$.
\end{theorem}

We are now ready to prove the Integral Basis Theorem.

\begin{proof}[Proof of \Cref{integral-basis}]
We prove this by induction on $|V|$. Let $P:=P(G)$.

\paragraph{Base case.} Suppose $G$ is BvN. We claim that $P$ has the integer decomposition property. To this end, let $x\in kP\cap \cZ^E$ for some integer $k\geq 1$, that is, $x\in \cZ^E_{\geq 0}$ and $x(\delta(v))=k,\, \forall v\in V$. If $k=1$, then $x$ is the incidence vector of a perfect matching. Otherwise, as $\frac{1}{k}x\in P$ and $P$ is an integral polytope, we can express $\frac{1}{k}x$ as a convex combination of the vertices of $P$, each of which is an incidence vector of a perfect matching of $G$. In particular, there is a perfect matching $M$ such that $\1_M\leq x$, where $\1_M\in \{0,1\}^E$ is the incidence vector of $M\subseteq E$. Let $x':=x-\1_M$, which satisfies $x'\in \cZ^E_{\geq 0}$ and $x'(\delta(v))=k-1,\, \forall v\in V$. By repeating this argument, we obtain a description of $x$ as the sum of $k$ vectors, each of which is an incidence vector of a perfect matching of $G$. Subsequently, $P$ has the integer decomposition property. As $G$ is matching-covered, the affine hull of $P$ is $\{x:x(\delta(v))=1,\forall v\in V\}$. Given that the GCD of the right-hand side values is $1$, it follows from \Cref{IDP-integral-basis} that $P\cap \cZ^E$ contains an integral basis for $\lin(P)$, as required.

\paragraph{Tight cuts.} Suppose there is a tight cut $C$. Let $G_1,G_2$ be the $C$-contractions of $G$, both of which are Petersen-free matching-covered graphs. By the induction hypothesis, for each $i\in [2]$, $\lin(P(G_i))$ has an integral basis $B^i$ that consists solely of the incidence vectors of some perfect matchings of $G_i$. Let $B:=B^1 \odot B^2$ which is a subset of $\cZ^E$ consisting of incidence vectors of some perfect matchings of $G$ that intersect $C$ just once, and is an integral basis for $\lin(P(G;C))$, by \Cref{composition}, part (3). Given that $C$ is a tight cut, it follows that $P(G;C)=P$, thus completing the induction step.

\paragraph{Non-BvN brick.} Otherwise, $G$ is a non-BvN non-Petersen brick. By the Intersection Theorem, there exist a separating facet-defining cut $C=\delta(X)$ and a perfect matching $M$ such that $|C\cap M|=3$. Let $G_1:=G/X$ and $G_2:=G/\bar{X}$, both of which are near-bricks, by \Cref{brick-count-near-brick}, part (3).
\medskip

Suppose in the first case that both $G_1,G_2$ are Petersen-free. Then by the induction hypothesis, for each $i\in [2]$, $\lin(P(G_i))$ has an integral basis $B^i$ that consists solely of the incidence vectors of some perfect matchings of $G_i$. Let $B':=B^1 \odot B^2$ which is a subset of $\cZ^E$ consisting of incidence vectors of some perfect matchings of $G$ that intersect $C$ just once, and is an integral basis for $\lin(P(G;C))$, by \Cref{composition}, part (3). We claim that $B:=B'\cup \{\1_M\}$ is an integral basis for $\lin(P)$. Linear independence is clear as $|M\cap C|>1$ while $b(C)=1$ for all $b\in B'$. Given that $\dim(P)=\dim(P(G;C))+1$, it therefore follows that $B$ is a basis for $\lin(P)$. It remains to prove that $B$ is an \emph{integral} basis. To this end, let $z\in \lin(P) \cap \cZ^E$, and express this vector as a unique linear combination of $B$: $z = \sum_{b\in B} \alpha_b b$. We need to show that $\alpha_b\in \cZ$ for all $b\in B$.

First we show that $\alpha_M:=\alpha_{\1_M}$ is an integer. As $z\in \lin(P)$, it follows that $z(\delta(u)) = z(\delta(v))$ for all $u,v\in V$. Let $c:=z(\delta(v))\in \cZ$. Recall that $C=\delta(X)$ for an odd-sized $X\subset V$. We have $$
|X| \cdot c = \sum_{v\in X} z(\delta(v)) = z(C) + 2z(E[X]),
$$ so $c \equiv z(C) \pmod{2}$. Given $v\in V$, we have $
c = z(\delta(v))= \sum_{b\in B} \alpha_b b(\delta(v)) = \sum_{b\in B} \alpha_b$, as each $b\in B$ is the incidence vector of a perfect matching. Thus, as $b(C)=1$ for all $b\in B'$ and $|M\cap C|=3$, $$z(C) = \sum_{b\in B} \alpha_b b(C) = 2\alpha_M + \sum_{b\in B} \alpha_b = 2\alpha_M + c.$$ Hence, $2\alpha_M = z(C)-c$, which is an even integer. Subsequently, $\alpha_M\in \cZ$. 

Thus, $z - \alpha_M \1_M$ is an integral vector in $\lin(B')$, so given that $B'$ is an integral basis for its linear hull, it follows that $\alpha_b\in \cZ$ for all $b\in B'$, as desired.
\medskip

Suppose in the remaining case that (at least) one of $G_1,G_2$, say $G_2=G/\bar{X}$, has a Petersen brick. We will adjust our choice of $C$ so that we fall in the previous case. 

We may assume that the near-brick $G_2$ is a Petersen brick. If not, then choose a minimal subset $Z\subset X$ such that $\delta(Z)$ is a tight cut of $G/\bar{X}$. In particular, $P(G;C)\subseteq P(G;\delta(Z))$. As $G$ is a brick, $\delta(Z)$ is not tight in $G$, so $\delta(Z)$ and $C$ define the same facet of $P$, so they are equivalent cuts and $G/Z/\bar{X}$ is bipartite, by \Cref{equivalent}. This, together with our minimal choice of $Z$, implies that $G/\bar{Z}$ is a brick. Subsequently, by changing $C=\delta(X)$ to an equivalent cut if necessary, we may assume $G_2$ is a Petersen brick.

Let $Y\subset X$ be the vertex set of a $5$-cycle of $G_2$. We claim that $\delta(Y)$ is a separating facet-defining cut of $G$ whose cut-contractions are Petersen-free, \b{and for which there exists a perfect matching $N$ such that $|N\cap \delta(Y)|=3$}.

\b{Observe that $|M\cap \delta(Y)|\in \{1,3\}$, and that the perfect matching $M$ can be redefined inside $G[X]$ to have either intersection size with $\delta(Y)$, whilst maintaining $|M\cap \delta(X)|=3$. Thus, there exists a perfect matching $N$ such that $|N\cap \delta(Y)|=3$}

Given that $\delta(Y)$ is separating in $G_2$ and $\delta(X)$ is separating in $G$, it follows that $\delta(Y)$ is separating in $G$.

To argue that $\delta(Y)$ is facet-defining for $G$, \rb{note first}{we may assume} that the perfect matching $M$ can be redefined inside $G[X]$, if necessary, such that $|M\cap \delta(Y)|=1$. Secondly, by \Cref{composition}, part (4) and \Cref{linear-basis-pete}, there exist perfect matchings $M_1,\ldots,M_d$ whose incidence vectors belong to $P(G;C)$ and form a basis for $\lin(P(G;C))$, where $|M_1\cap \delta(Y)|=5$, $|M_i\cap \delta(Y)|=1$ for $i\in \{2,\ldots,d\}$, and every edge of $G$ belongs to one of $M_2,\ldots,M_d$. By swapping $M_1$ with $M$, we obtain $d$ perfect matchings $M,M_2,\ldots,M_d$ whose incidence vectors belong to $P(G;\delta(Y))$ and are linearly independent, because $|M\cap C|=3$ while $|M_i\cap C|=1$ for $i\in \{2,\ldots,d\}$. Thus, $P(G;\delta(Y))$ is a facet of $P$.

Finally, \b{after replacing $C$ by $\delta(Y)$,} the near-brick $G_i$ must be Petersen-free because it contains a triangle, and such a triangle also belongs to the unique brick of $G_i$, for each $i\in [2]$.
\end{proof}

We are now ready to prove the Main Theorem.

\begin{proof}[Proof of \Cref{lattice-basis}]
We proceed by induction on $|V|$. For the base case, suppose $G$ is a brick. If $G$ is a non-Petersen brick, then the result follows from \Cref{integral-basis}. Otherwise, $G$ is a Petersen brick, and the result can be readily checked. 

For the induction step, suppose $C$ is a tight cut of $G$, and let $G_1,G_2$ be the $C$-contractions, where $G_1$ has Petersen bricks $H_1,\ldots,H_{q}$, and $G_2$ has Petersen bricks $H_{q+1},\ldots,H_p$, obtained through tight cut decompositions of $G_1,G_2$, respectively. Let $L_i:=L(G_i)$ and $\bar{L}_i:=\lin(P(G_i))\cap \cZ^{E(G_i)}$, for $i\in [2]$. By the induction hypothesis, $L_i$ has a lattice basis $B_i$ consisting solely of some perfect matchings of $G_i$, for $i=1,2$. Furthermore, 
\begin{align*}
L_1 &= \bar{L}_1 \cap \left\{x : x(A_i) \equiv 0 \pmod{2},~i=1,\ldots,q\right\},\\
L_2 &= \bar{L}_2 \cap \left\{y : y(A_i) \equiv 0 \pmod{2},~i=q+1,\ldots,p\right\},
\end{align*}
where each $A_i$ is the edge set of some $5$-cycle of $H_i$. Let $B:=B_1\odot B_2$ which, \b{by \Cref{composition}, part (1)}, consists of some perfect matchings of $G$. By \Cref{composition}, part (2), $B$ is a lattice basis for $L$. Furthermore, we claim that the above set equalities imply $$
L = \bar{L} \cap \left\{z : z(A_i) \equiv 0 \pmod{2},~\forall i\in [p]\right\}.
$$ $(\subseteq)$ is clear. $(\supseteq)$: Pick $f\in \bar{L}$ such that $f(A_i)\equiv 0\pmod{2},~\forall i\in [p]$. Then $f=x\odot y$ for $x\in \bar{L}_1$ and $y\in \bar{L}_2$, where clearly $x(A_i) \equiv 0 \pmod{2},~i=1,\ldots,q$ and $y(A_i) \equiv 0 \pmod{2},~i=q+1,\ldots,p$. Thus, $x\in L_1,y\in L_2$, so $f$ belongs to $L$, by \Cref{composition}, part (2). This completes the induction step.
\end{proof}

\section*{Acknowledgements and declarations}

Parts of this work were carried out when AA was a Visiting Professor hosted generously by the G-SCOP Laboratory, Université Grenoble Alpes in Fall 2024, and when OS was visiting the London School of Economics in Winter 2025. This work was supported in part by EPSRC grant EP/X030989/1. We would like to thank G\'{e}rard Cornu\'{e}jols and Andr\'{a}s Seb\H{o} for fruitful discussions about this work.
\b{Finally, we thank a reviewer whose comments improved the presentation of the paper.}

\paragraph{Data availability statement.} No data are associated with this article. Data sharing is not applicable to this article.

\paragraph{Conflicts of interests/competing interests statement.} The authors have no conflicts of interest nor competing interests to declare that are relevant to the content of this article.

{\small 
\bibliographystyle{alpha}
\bibliography{Petersen-graph-references}}

\end{document}